\documentclass[12pt]{amsart} 
\usepackage{lipsum}
\makeatletter
\g@addto@macro{\endabstract}{\@setabstract}
\newcommand{\authorfootnotes}{\renewcommand\thefootnote{\@fnsymbol\c@footnote}}%
\makeatother

\usepackage{enumitem}
\usepackage[dvipsnames]{xcolor}
\usepackage{setspace}
\usepackage[english]{babel}
\usepackage{amsmath,amsfonts,amssymb,amsthm,epsfig,epstopdf,titling,url,array}
\usepackage[margin=1.2in]{geometry}
\usepackage{tikz-cd}
\usepackage{calc}
\usepackage{color,hyperref}
\usepackage[noabbrev]{cleveref}
\hypersetup{colorlinks}
\hypersetup{colorlinks={true},linkcolor={red},citecolor=green}
\usepackage{graphicx}
\usepackage{mathrsfs}

\graphicspath{{images/}}
\theoremstyle{plain}
\newtheorem{theorem}{Theorem}

\newtheorem{thm}{Theorem}[subsection]
\newtheorem{lem}[thm]{Lemma}

\newcommand\scalemath[2]{\scalebox{#1}{\mbox{\ensuremath{\displaystyle #2}}}}
\theoremstyle{plain}
\newtheorem{defn}[thm]{Definition}

\theoremstyle{remark}
\newtheorem{rem}[thm]{Remark}

\newcommand{\mycomment}[1]{}
\usepackage{tikz}

\makeatletter
\usepackage[new]{old-arrows}
\usepackage{bm} 
\usepackage{multicol}
\usetikzlibrary{decorations.pathmorphing}
\DeclareMathOperator{\GL}{GL}

\DeclareMathOperator{\ch}{ch}

\usepackage{mathtools}

\newcommand\blfootnote[1]{%
  \begingroup
  \renewcommand\thefootnote{}\footnote{#1}%
  \addtocounter{footnote}{-1}%
  \endgroup
}

\makeatletter
\renewcommand{\tocsection}[3]{%
  \indentlabel{\@ifnotempty{#2}{\bfseries\ignorespaces#1 #2\quad}}\bfseries#3}
\renewcommand{\tocsubsection}[3]{%
  \indentlabel{\@ifnotempty{#2}{\ignorespaces#1 #2\quad}}#3}

\renewcommand{\tocsubsubsection}[3]{%
  \indentlabel{\@ifnotempty{#2}{\ignorespaces#1 #2\quad}}#3}

\newcommand\@dotsep{4.5}
\def\@tocline#1#2#3#4#5#6#7{\relax
  \ifnum #1>\c@tocdepth 
  \else
    \par \addpenalty\@secpenalty\addvspace{#2}%
    \begingroup \hyphenpenalty\@M
    \@ifempty{#4}{%
      \@tempdima\csname r@tocindent\number#1\endcsname\relax
    }{%
      \@tempdima#4\relax
    }%
    \parindent\z@ \leftskip#3\relax \advance\leftskip\@tempdima\relax
    \rightskip\@pnumwidth plus1em \parfillskip-\@pnumwidth
    #5\leavevmode\hskip-\@tempdima{#6}\nobreak
    \leaders\hbox{$\m@th\mkern \@dotsep mu\hbox{.}\mkern \@dotsep mu$}\hfill
    \nobreak
    \hbox to\@pnumwidth{\@tocpagenum{\ifnum#1=1\bfseries\fi#7}}\par
    \nobreak
    \endgroup
  \fi}
\AtBeginDocument{%
\expandafter\renewcommand\csname r@tocindent0\endcsname{0pt}
}
\def\l@subsection{\@tocline{2}{0pt}{2.5pc}{5pc}{}}
\def\l@subsubsection{\@tocline{2}{0pt}{4pc}{5pc}{}}
\makeatother

\date{} 

\begin{document}
\hypersetup{citecolor=blue}
\hypersetup{linkcolor=red}
\begin{center}
  \normalsize
A NOTE ON TORIC PERIODS IN UNRAMIFIED FAMILIES
  \normalsize
  \\
  \bigskip
  Alexandros Groutides \par 
  
Instituto Matem\'aticas de la Universidad de Sevilla\par
\end{center}
\begin{abstract}
    Let $A$ be the algebra $\mathbf{C}[X_1^{\pm 1},X_2^{\pm 1}]$ and $Q(A)$ its quotient field. In this short article we exhibit the correct normalization for the toric period on the parabolically induced unramified family over $Q(A)$, so that it behaves optimally under restriction to the parabolically induced unramified family over $A$. Finally, we explicitly describe the structure of the $A$-module given by the image of this normalized period. These results answer questions raised by D. Prasad, and point towards potential generalizations to a broader unramified Gan-Gross-Prasad setting.
\end{abstract}
     \blfootnote{We gratefully acknowledge support from the following research grants: ERC 101001051: Shimura varieties and the Birch–Swinnerton-Dyer conjecture, ATR 2024-154613: La correspondencia de
Langlands local y el dual unitario de los grupos $p$-ádicos clásicos.}
\vspace{-2em}

\section{Introduction}
\subsection{Setup}Let $p$ be a prime and $F$ a finite field extension of $\mathbf{Q}_p$ with ring of integers $\mathcal{O}$. Write $\varpi$ for a fixed choice of uniformizer and $q$ for the cardinality of the residue field of $F$. Let $A:=\mathbf{C}[X_1^{\pm 1},X_2^{\pm 1}]$ be the commutative $\mathbf{C}$-algebra of complex Laurent polynomials in the variables $X_1$ and $X_2$, and $Q(A)$ its quotient field. Let $T$ denote the diagonal torus of $\GL_2$, $B$ the upper triangular Borel subgroup with modular character $\delta_B$, and $N$ the corresponding unipotent radical. 

Let $\chi:T(F)\rightarrow A$ be the unramified character sending the matrix $\mathrm{diag}(\varpi^a,\varpi^b)$ to $X_1^aX_2^b$, which we regard as a character of $B(F)$ in the usual manner. We will consider the two normalized parabolically induced families, one defined over $Q(A)$ and the other over $A$:
$$
    \scalemath{0.9}{i_{Q(A)}(\chi):=\left\{f:\mathrm{GL}_2(F)\rightarrow Q(A)\ \mathrm{smooth,}\ \mathrm{and}\ f(bg)=(\delta_B^{1/2}\chi)(b)f(g)\ \forall\  b\in B(F),g\in\mathrm{GL}_2(F)\right\} }$$
    $$ \scalemath{0.9}{
    i_{A}(\chi):=\left\{f:\mathrm{GL}_2(F)\rightarrow A\ \mathrm{smooth,}\ \mathrm{and}\ f(bg)=(\delta_B^{1/2}\chi)(b)f(g)\ \forall\  b\in B(F),g\in\mathrm{GL}_2(F)\right\}.}
$$
The $Q(A)[\GL_2(F)]$-module $i_{Q(A)}(\chi)$ is irreducible, whereas the $A[\GL_2(F)]$-module $i_A(\chi)$ is not. We write $f^\mathrm{sph}$ for the normalized $\GL_2(\mathcal{O})$-invariant vector defined via the Iwasawa decomposition, as $f^\mathrm{sph}(bk):=(\chi\delta_B^{1/2})(b)$ for all $b\in B(F)$, $k\in \GL_2(\mathcal{O})$. 
\subsection{Motivation and main result} Under the usual embedding $\GL_1\hookrightarrow \GL_2,\ a\mapsto\mathrm{diag}(a,1)$, the space of unramified toric periods
\begin{align}\label{eq: hom space}\mathrm{Hom}_{\GL_1(F)}(i_{Q(A)}(\chi),Q(A))\end{align}
is one-dimensional over $Q(A)$, and any such non-zero linear form is non-vanishing on $f^\mathrm{sph}$. This is due to celebrated results of Waldspurger \cite{waldspurger1985valeurs} and Prasad \cite{prasad1990trilinear} over $\mathbf{C}$, which remain valid over the field $Q(A)$ which contains $\mathbf{C}$. To summarize, there is a unique up to $Q(A)$-scalars linear form in \eqref{eq: hom space}, and it is non-vanishing on $f^\mathrm{sph}$.

Now, since $i_A(\chi)$ is defined over $A$, one would like to show that there exists an explicit normalized toric period $l_A$ in \eqref{eq: hom space}, unique up to scalar multiplication by $A^\times$, for which $l_A(i_A(\chi))$ is contained in $A$ and $a^{-1}l_A(i_A(\chi))$ is not contained in $A$ for any non-zero $a\notin A^\times$. {In other words, $l_A(i_A(\chi))$ is contained in $A$ but is not contained in any proper principal ideal of $A$.} Having identified $l_A$ it is then natural to enquire about the precise structure of the $A$-module $l_A(i_A(\chi))\subseteq A$. These are the problems we address. 

\begin{theorem}\label{thm intro main} Let $l_A$ be the unique linear form in \eqref{eq: hom space} normalized to satisfy $l_A(f^\mathrm{sph})=1-q^{-1}X_1X_2^{-1}.$
Then, the $A$-module $l_A(i_A(\chi))$ is precisely the non-free $A$-submodule of $A$ given by
    $$(1-q^{-1/2}X_1)A+(1-q^{-1}X_1X_2^{-1})A=(1-q^{1/2}X_2)A+ (1-q^{-1}X_1X_2^{-1})A\subseteq A.$$
    In particular, $l_A(i_A(\chi))$ is contained in $A$ and $a^{-1}l_A(i_A(\chi))$ is not contained in $A$ for any non-zero $a\notin A^\times$.
\end{theorem}
\subsection{Some remarks}
\begin{rem}
    Let us fix the standard basis $\{e_1,e_2\}$ for $X^*(T)$, and write $\alpha:=e_1-e_2$ for the positive root. Then the normalized value of $l_A(f^\mathrm{sph})$ in \Cref{thm intro main} which is given by
    \begin{align*}1-q^{-1}X_1X_2^{-1}=\prod_{\alpha>0}\left(1-q^{-1}\chi(\alpha^\vee(\varpi))\right),\end{align*}
    is precisely the normalization factor appearing in the Casselman-Shalika formula \cite{casselman1980unramifiedII}. With this in mind, we expect analogs of the integrality theorem above to hold in much more general Gan-Gross-Prasad settings for unramified representations of $G\times H$, admitting non-zero $H$-invariant linear forms valued in $\mathbf{C}(X_1,\dots, X_r)$ where $r$ denotes the split rank of $G\times H$. We hope that this short article will spark interest for further work in this direction which is largely unexplored. 
    
    {In the present article, we study the case $G\times H=\GL_2(F)\times \GL_1(F)$, where the $G$-representation is $i_{Q(A)}(\chi)$, the $H$-representation is the trivial one, and $H$ is regarded as a subgroup of $G\times H$ via the diagonal embedding. In this case $r=3$, but since the $H$-representation is trivial, the representations and the linear forms considered are defined over $Q(A)=\mathbf{C}(X_1,X_2)\subseteq \mathbf{C}(X_1,X_2,X_3)$. Our techniques can easily be applied to include the slightly more general case where the $H$-representation is the unramified family $F^\times=\GL_1(F)\rightarrow \mathbf{C}(X_3)$ defined by $\varpi\mapsto X_3.$}
\end{rem}
\begin{rem}
    We also point out an interesting phenomenon arising from the fact that $i_A(\chi)$ fails to be irreducible as an $A[\GL_2(F)]$-module. Let $l_1$ denote the period in \eqref{eq: hom space} such that $l_1(f^\mathrm{sph})=1$. Then, for example using \cite[Theorem $3.1.1$]{groutides2024integral}, and without needing any knowledge of what $l_1$ looks like, other than the fact that it maps $f^\mathrm{sph}$ to $1$, we see that $l_1(A[\GL_2(F)]\cdot f^\mathrm{sph})$ is contained in $A$. However, \Cref{thm intro main} clearly implies that the whole of $l_1(i_A(\chi))$ (which contains $l_1(A[\GL_2(F)]\cdot f^\mathrm{sph})$) is not itself contained in $A$. We also expect this to be the case for higher rank groups.

    Finally, it follows easily from the defining property of the linear form $l_A$ and the multiplicity-one over $Q(A)$, that the space $\mathrm{Hom}_{A[\GL_1(F)]}(i_A(\chi),A)$ is free of rank one over $A$, with $l_A$ as a basis. So, unlike the situation over $Q(A)$ where any basis element of $\mathrm{Hom}_{Q(A)[\GL_1(F)]}(i_{Q(A)}(\chi),Q(A))\simeq Q(A)$ has image exactly $Q(A)$, we see that any basis element of $\mathrm{Hom}_{A[\GL_1(F)]}(i_A(\chi),A)\simeq A$ has non-free image over $A$. Once again we expect this to be the case in higher ranks. 
\end{rem}

\subsection{Acknowledgements} I would like to thank Dipendra Prasad for bringing this problem to my attention and for valuable discussions. I would also like to thank Justin Trias for many insightful comments.
{
  \hypersetup{linkcolor=black}
  \setcounter{tocdepth}{2}
  \tableofcontents
}
\section{Realizing linear forms}\label{sec: realizing ell}
\subsection{The Whittaker functional}
For $i\in\mathbf{Z}$, we define the open compact subgroup $N_i(\mathcal{O}):=\left[\begin{smallmatrix}
    1 & \varpi^i\mathcal{O}\\
    &1
\end{smallmatrix}\right]\subseteq N(F)$. 
Let $\psi:F\rightarrow\mathbf{C}^\times$ be {an} additive character of conductor $\mathcal{O},$ regarded as a $Q(A)$-valued character of the unipotent radical $N(F)$ in the obvious way, and let $w{:=\left[\begin{smallmatrix}
    & 1 \\
    1& 
\end{smallmatrix}\right]}$ be the non-trivial Weyl element of $\GL_2$.  We want to realize a Whittaker functional in $\mathrm{Hom}_{N(F)}(i_{Q(A)}(\chi),\psi)$ as Casselman-Shalika do in \cite{casselman1980unramifiedII}. Let $I_1$ denote the $B(F)$-stable subspace of functions with support contained in the big Bruhat cell $B(F)wB(F)$. We then write $\Lambda_\chi$ for the functional defined on $I_1$ as
$$\Lambda_\chi(f):=\int_{N(F)}f(wn)\psi^{-1}(n)\ dn,\ f\in I_1.$$
As in \textit{op.cit.}, for any $f\in I_1$, the function $n\mapsto f(wn)$ belongs {to} $C_c^\infty(N(F),Q(A))${,} the space of locally constant and compactly supported $Q(A)$-valued functions on $N(F)$. {Here, the integral $\int_{N(F)} dn$} denotes the unique translation-invariant linear {functional} $C_c^\infty(N(F),Q(A))\rightarrow Q(A)$ which sends the characteristic function $\ch_{N(\mathcal{O})}$ to $1.$ 

Having defined $\Lambda_\chi$ on $I_1$, we can extend it to all of $i_{Q(A)}(\chi)$  using exactness of the $\psi$-twisted Jacquet functor which holds over $Q(A),$ and the exact sequence of $Q(A)[B(F)]$-modules
\begin{align}\label{eq: exact seq}0\rightarrow I_1\rightarrow i_{Q(A)}(\chi)\rightarrow
\delta_B^{1/2}\chi \rightarrow 0
\end{align}
where the surjection is given by restriction of functions to $B(F).$  Thus, we can now identify $i_{Q(A)}(\chi)$ with its Whittaker model, via 
\begin{align}\label{eq: Whittaker}i_{Q(A)}(\chi)\simeq \mathcal{W}(i_{Q(A)}(\chi),\psi),\ f\mapsto (W_f:g\mapsto \Lambda_\chi(g\cdot f)).\end{align}
\subsection{The formal zeta-integral}
Now, we follow the approach of \cite[\S $3.1$]{kurinczuk2017rankin} to define a linear form on this Whittaker model $\mathcal{W}(i_{Q(A)}(\chi),\psi)$ using a formal algebraic version of the classical $\GL_2$ zeta-integral over $\mathbf{C}.$
\begin{defn}\label{def: linear form whitt}
Let $\lambda:\mathcal{W}(i_{Q(A)}(\chi),\psi)\rightarrow Q(A)$ be given by 
 $$\lambda(W):=\frac{I(W,X)}{L(i_{Q(A)}(\chi),X)}\biggr\rvert_{X=q^{-1/2}}$$
 where $I(W,X)\in Q(A)((X))$ is the formal Laurent series in the variable $X$
 $$I(W,X):=\sum_{k\in\mathbf{Z}}\left\{\int_{\mathcal{O}^\times}W(\left[\begin{smallmatrix}
     a\varpi^k & \\
     & 1
 \end{smallmatrix}\right])\ d^\times a \right\} q^{k/2}X^k,$$
 $L(i_{Q(A)}(\chi),X)$ is the $L$-factor $\frac{1}{(1-X_1X)(1-X_2X)}$, and the translation-invariant linear functional $\int_{\mathcal{O}^\times}d^\times a$ on $C_c^\infty(\mathcal{O}^\times,Q(A))$ is normalized to send $\ch_{\mathcal{O}^\times}$ to $1.$
\end{defn}
This is well-defined since the Whittaker coefficients in the curly brackets are zero for $k\ll 0$, and the $I(W,X)$ span the fractional ideal $ L(i_{Q(A)}(\chi),X)\cdot Q(A)[X,X^{-1}]$ as we vary $W\in \mathcal{W}(i_{Q(A)}(\chi),\psi)$. These either follow from the general arguments in \cite[\S $3.1$]{kurinczuk2017rankin}, or in this spherical $\GL_2$ case, from the Shintani formula for the spherical Whittaker function \cite{shintani1976explicit} which remains valid over $Q(A).$ It now follows from an easy direct computation that the linear form $\lambda$ is in fact an element of $$\mathrm{Hom}_{\GL_1(F)}(\mathcal{W}(i_{Q(A)}(\chi),\psi),Q(A)).$$ 
Moreover, as in the complex theory, if we write $W^\mathrm{sph}$ for the normalized spherical vector in $\mathcal{W}(i_{Q(A)}(\chi),\psi)$ satisfying $W^\mathrm{sph}(1)=1$, we have {$I(W^\mathrm{sph},X)=L(i_{Q(A)}(\chi),X)$ and hence} $\lambda(W^\mathrm{sph})=1.$ We can now descend back to the principal-series model and realize $l_A$ as follows.
 \begin{defn}
     Let $l_A$ be the linear form $$l_A:i_{Q(A)}(\chi)\rightarrow Q(A),\ f\mapsto \lambda(W_f)$$
     where $\lambda$ is the linear form of \Cref{def: linear form whitt} and $W_f$ is the Whittaker function associated to $f$ as in \eqref{eq: Whittaker}.
 \end{defn}

\begin{lem}
   We have $l_A(f^\mathrm{sph})=1-q^{-1}X_1X_2^{-1}.$
    \begin{proof}
       This follows from the fact that $W_{f^\mathrm{sph}}=(1-q^{-1}X_1X_2^{-1})W^\mathrm{sph}$, which is a well-known fact over $\mathbf{C}$. Over $Q(A)$, this equality remains true. Indeed, if we use the algebraic approach of \cite{casselman1980unramifiedII} which works verbatim over $Q(A)$, and let $\mathcal{P}_{f^\mathrm{sph}}:\GL_2(F)\rightarrow Q(A)$ denote the function
       $$\mathcal{P}_{f^\mathrm{sph}}(g):=q^{-1}\int_{N_{-1}(\mathcal{O})}\psi^{-1}(n)f^\mathrm{sph}(gn)\ dn,$$
       then it is clear that it has support contained inside the big cell $B(F)wB(F)$. More precisely, it has support contained in $B(F)w N_{-1}(\mathcal{O}).$ Combining this with \eqref{eq: exact seq}, we have $$\Lambda_\chi(f^\mathrm{sph})=\Lambda_\chi(\mathcal{P}_{f^\mathrm{sph}})=\int_{N_{-1}(\mathcal{O})}f^\mathrm{sph}(wn)\psi^{-1}(n)\ dn=1-q^{-1}X_1X_2^{-1}.$$
    \end{proof}
\end{lem}
\section{Restricting to $i_A(\chi)$}
\subsection{The inclusion} The goal of this section is to prove the following result.
\begin{thm}\label{thm 3.1.1}
    Let $l_A$ be the unique linear form in $\mathrm{Hom}_{\GL_1(F)}(i_{Q(A)}(\chi),Q(A))$ normalized to send $f^\mathrm{sph}$ to $1-q^{-1}X_1X_2^{-1}$. Then, the $A$-module $l_A(i_A(\chi))$ is contained in the  non-free $A$-submodule of $A$ given by 
    \begin{align}\label{eq: A-mod}(1-q^{-1/2}X_1)A+(1-q^{-1}X_1X_2^{-1})A=(1-q^{1/2}X_2)A+ (1-q^{-1}X_1X_2^{-1})A\subseteq A.\end{align}
\begin{proof}
Let $f:\GL_2(F)\rightarrow A$ be an element of $i_A(\chi)$ and set $a_f:=f(1)\in A$. It follows from the Bruhat decomposition $$\GL_2(F)=B(F)\sqcup B(F)wB(F)$$
that $f_w:=f-a_ff^\mathrm{sph}\in i_A(\chi)$ has support contained in the big cell $B(F)wB(F)$. Hence, for any $a\in F^\times$, the function $\left[\begin{smallmatrix}
    a & \\
    & 1
\end{smallmatrix}\right]\cdot f_w$ also has support on the big cell. Thus, from the defining property of the associated Whittaker function, we have 
\begin{align}\label{eq: 2}
    I(&W_{f_w},X)=\sum_{k\in\mathbf{Z}}\left\{\int_{\mathcal{O}^\times}W_{f_w}(\left[\begin{smallmatrix}
    a\varpi^k & \\
    & 1
\end{smallmatrix}\right])\ d^\times a\right\}q^{k/2}X^k\\
\nonumber&=\sum_{k\in\mathbf{Z}}\left\{\int_{\mathcal{O}^\times}\int_{N(F)}f_w(wn\left[\begin{smallmatrix}
    a\varpi^k & \\
    & 1
\end{smallmatrix}\right])\psi^{-1}(n)\ dn\ d^\times a\right\}q^{k/2}X^k\\
\nonumber&=\sum_{k\in\mathbf{Z}}(\chi\delta_B^{1/2})(\left[\begin{smallmatrix}
    1 & \\
    & \varpi^k
\end{smallmatrix}\right])\left\{\int_{\mathcal{O}^\times}\int_{N(F)}f_w(w\left[\begin{smallmatrix}
    a\varpi^{k} & \\
    & 1
\end{smallmatrix}\right]^{-1}n\left[\begin{smallmatrix}
    a\varpi^k & \\
    & 1
\end{smallmatrix}\right])\psi^{-1}(n)\ dn\ d^\times a\right\}q^{k/2}X^k\\
\nonumber&=\sum_{k\in\mathbf{Z}} \left\{\int_{\mathcal{O}^\times}\int_{N(F)}f_w(wn)\psi_{a\varpi^k}^{-1}(n)\ dn\ d^\times a\right\} X_2^kX^k
\end{align}
where $\psi_{a\varpi^k}(n):=\psi(\left[\begin{smallmatrix}
    a\varpi^k & \\
    & 1
\end{smallmatrix}\right]n\left[\begin{smallmatrix}
    a\varpi^k & \\
    & 1
\end{smallmatrix}\right]^{-1})$. Note that the last equality follows from the fact that the linear functional $\int_{N(F)} dn$ gets scaled by $|a\varpi^k|=q^{-k}$ under conjugation by $\left[\begin{smallmatrix}
    a\varpi^k & \\
    & 1
\end{smallmatrix}\right]^{-1}.$ Since $n\mapsto f_w(wn)$ is smooth and compactly supported, there exists an integer $c$ only depending on $f_w$, such that $$\int_{N(F)}f_w(wn)\psi_{a\varpi^k}^{-1}(n)\ dn=\sum_{i=1}^mf_w(wn_i)\int_{N_c(\mathcal{O})}\psi_{a\varpi^k}^{-1}(n_in)\ dn$$for some $n_i\in N(F)$. By assumption, the coefficients $f_w(wn_i)$ are contained in $A$. Thus, \eqref{eq: 2} is a finite $A$-linear combination of expressions of the form 
\begin{align}\label{eq: 3}
    \sum_{k\in\mathbf{Z}}\left\{\int_{\mathcal{O}^\times}\psi_{a\varpi^k}^{-1}(n_0)\int_{N_c(\mathcal{O})}\psi_{a\varpi^k}^{-1}(n)\ dn\ d^\times a\right\} X_2^kX^k
\end{align}
for some $c\in\mathbf{Z}$ and some $n_0\in N(F).$ But, for $k\ll 0$, $\psi_{a\varpi^k}(n)$ is a non-trivial character of the compact group $N_c(\mathcal{O})$ and hence the inner integral vanishes. Alternatively, if the inner integral does not vanish, it is identically equal to the fixed quantity $\int_{N_c(\mathcal{O})}1\ dn\in\mathbf{Q}$. Moreover, for $k\gg 0$, $\psi_{a\varpi^k}^{-1}(n_0)=1$ for all $a\in\mathcal{O}^\times$. It now follows that \eqref{eq: 3}, and hence \eqref{eq: 2}, is an element of $\frac{1}{1-X_2X}A.$ Thus, 
\begin{align*}
    l_A(f_w)=\lambda(W_{f_w})\in (1-q^{-1/2}X_1)A,
\end{align*}
and hence 
\begin{align*}  l_A(f)&=a_fl_A(f^\mathrm{sph})+l_A(f_w)\\
&=a_f(1-q^{-1}X_1X_2^{-1})+l_A(f_w)\in (1-q^{-1}X_1X_2^{-1})A+(1-q^{-1/2}X_1)A\subseteq A.
\end{align*}
\end{proof}
\end{thm}
\subsection{The equality}
 We now strengthen \Cref{thm 3.1.1} by showing that the $A$-module $l_A(i_A(\chi))$ is actually the whole of \eqref{eq: A-mod}. In other words, we prove the following result.
 \begin{thm}\label{thm main}
    Let $l_A$ be the unique linear form in $\mathrm{Hom}_{\GL_1(F)}(i_{Q(A)}(\chi),Q(A))$ normalized to send $f^\mathrm{sph}$ to $1-q^{-1}X_1X_2^{-1}$. Then, the $A$-module $l_A(i_A(\chi))$ is precisely the non-free $A$-submodule of $A$, given by
    $$(1-q^{-1/2}X_1)A+(1-q^{-1}X_1X_2^{-1})A=(1-q^{1/2}X_2)A+ (1-q^{-1}X_1X_2^{-1})A\subseteq A.$$
    In particular, $l_A(i_A(\chi))$ is contained in $A$ and $a^{-1}l_A(i_A(\chi))$ is not contained in $A$ for any $a\notin A^\times$.
\begin{proof}
 Since $l_A(f^\mathrm{sph})=1-q^{-1}X_1X_2^{-1}$, it suffices to show that there exists some $f\in i_A(\chi)$ for which $l_A(f)=1-q^{-1/2}X_1.$ Let 
 $$I(\mathcal{O}):=\left[\begin{matrix}
     \mathcal{O} & \mathcal{O}\\
     \varpi\mathcal{O} & \mathcal{O}
 \end{matrix}\right]\cap \GL_2(\mathcal{O})$$
 denote the Iwahori subgroup of $\GL_2(F)$. Let $f_0$ denote the following element of $i_A(\chi)$
 \begin{align*}
     f_0(g):=\begin{dcases}
         0,\ &\mathrm{if}\ g\not\in B(F)wI(\mathcal{O})\\
         (\chi\delta_B^{1/2})(b),\ &\mathrm{if}\ g=bwi\in B(F)wI(\mathcal{O}).
     \end{dcases}
 \end{align*} Note that this is the function denoted by $\phi_w$ in \cite[\S $2$]{casselman1980unramified}. First, note that the support of $f_0$ is precisely $B(F)wI(\mathcal{O})\subseteq B(F)wB(F)$, where the inclusion follows from \cite[Proposition $1.3$]{casselman1980unramified}. Thus, for any $a\in F^\times$, the Whittaker function $W_{f_0}(\left[\begin{smallmatrix}
     a & \\
      & 1
 \end{smallmatrix}\right])$ is given by the unipotent integral of the previous section. Thus, as in \eqref{eq: 2} we have
\nonumber \begin{align}
     I(W_{f_0},X)&=\sum_{k\in\mathbf{Z}} \left\{ \int_{\mathcal{O}^\times}\int_{N(F)} f_0(wn)\psi_{a\varpi^k}^{-1}(n)\ dn\ d^\times a \right\} X_2^k X^k\\
     &=\sum_{k\in\mathbf{Z}} \left\{ \int_{\mathcal{O}^\times}\int_{N(\mathcal{O})}\psi_{a\varpi^k}^{-1}(n)\ dn\ d^\times a \right\} X_2^k X^k\\
     &=\frac{1}{1-X_2X}
 \end{align}
 where the last equality follows once again by orthogonality, and the second equality follows from the fact that for $n\in N(F)$, we have $wn\in B(F)wI(\mathcal{O})$ if and only if $n\in N(\mathcal{O})$. Now, we can conclude that 
 $$l_A(f_0)=1-q^{-1/2}X_1.$$
 as required.
 \end{proof}
 \end{thm}

\noindent\textit{Instituto Matemáticas de la Universidad de Sevilla, Seville 41012,
Spain}.\\
\textit{Email address}: alexgroutides21@gmail.com
\end{document}